\documentclass[11pt]{amsart}
 
\usepackage{amsfonts}
\usepackage{amssymb}
\usepackage{amscd}
\input xy
\xyoption{all}
 
\newcommand{\marginlabel}[1]%
  {\mbox{}\marginpar{\raggedleft\hspace{0pt}\bfseries\sf#1}}

\def\NN{{\mathbb N}}

\def\CC{{\mathbb C}}
\def\AA{{\mathbb A}}
\def\RR{{\mathbb R}}
\def\QQ{{\mathbb Q}}

\def\be{\mathbf{e}}

\def\be{\mathbf{e}}

\DeclareMathOperator{\codim}{codim}

\DeclareMathOperator{\ord}{ord}
\DeclareMathOperator{\mld}{mld}
\DeclareMathOperator{\Image}{Im}

\newtheorem{lemma}{Lemma}[section]
\newtheorem{theorem}[lemma]{Theorem}
\newtheorem{corollary}[lemma]{Corollary}

\theoremstyle{definition}

\newtheorem{remark}[lemma]{Remark}

\theoremstyle{remark}
\newtheorem*{remark*}{Remark}
\newtheorem*{note*}{Note}

\begin{document}

\title{Inversion of Adjunction for
local complete intersection varieties}

\author[L. Ein]{Lawrence Ein}
\address{Department of Mathematics, Statistics and
Computer Science, University of Illinois at Chicago,
851 Morgan St., M/C. 249, Chicago, IL 60607-7045, USA}
\email{{\tt ein@math.uic.edu}}

\author[M. Musta\c{t}\v{a}]{Mircea~Musta\c{t}\v{a}}
\address{Department of Mathematics, Harvard University,
One Oxford Street, Cambridge, MA 02138, USA}
\email{{\tt mirceamustata@yahoo.com}}

\thanks{\mbox {  } \mbox {  } 2000
\newline \mbox {  } \mbox {  }
{\it Key words and phrases}. Jet schemes, minimal log discrepancies,
inversion of adjunction.}

\begin{abstract}
We extend Inversion of Adjunction and Semicontinity of
minimal log discrepancies from the case of a smooth ambient variety,
as in \cite{EMY}, to the case of local
complete intersection varieties.
\end{abstract}

\maketitle

\section{Introduction}

In this note we build on our previous work with Takehiko Yasuda in \cite{EMY},
to prove a precise version of Inversion of Adjunction for varieties
which are locally complete intersections. 
\cite{EMY} gives a description of certain
invariants of singularities, known as
minimal log discrepancies, using the codimension of suitable spaces of arcs.
This was enough to give Inversion of Adjunction
in the case when the ambient variety is smooth. 
The main new ingredient in this paper
is the study of the equations defining the jets
which can be lifted to the space of arcs.

Recall the setting for Inversion of Adjunction.
Consider a pair $(X,Y)$, where $X$ is a $\QQ$-Gorenstein normal variety, and
where $Y$ stands for a formal combination
$\sum_{i=1}^kq_i\cdot Y_i$, where $q_i\in\RR_+$
and $Y_i\subset X$ are proper closed subschemes. 
Given such a pair and a closed subset $W\subseteq X$,
one can get an invariant $\mld(W;X,Y)$, called the minimal
log discrepancy. Loosely speaking, this invariant 
comes from evaluating $K_X+Y$ along divisors over $X$, with
centers in $W$. We refer to Section 3 for the precise definition.
The following is our main result.

\begin{theorem}\label{inv_adj0}
Let $X$ be a normal, local complete intersection variety, 
and $Y=\sum_i q_i\cdot Y_i$, where $q_i\in\RR_+$ and $Y_i\subset X$
are proper closed subschemes. If 
$D\subset X$ is
a normal effective  Cartier divisor such that 
$D\not\subseteq\bigcup_iY_i$, then
for every proper closed subset $W\subset D$, we have
$$\mld(W;X,D+Y)=\mld(W;D,Y\vert_D).$$
\end{theorem}

It is a conjecture due to Koll\'ar and Shokurov
that the equality in the above conjecture is true if $X$ is just
normal and $\QQ$-Gorenstein. Note that the equality ``$\leq$''
is easy and is known in general. We refer to \cite{kollar}
for a discussion of this conjecture, related topics, and for the
proof of some special cases via vanishing theorems.
Our approach will be via spaces of arcs, based on the results from
\cite{EMY}, where the theorem was proved when $X$ is smooth.

Ambro proposed in \cite{ambro} a conjecture on the
semicontinuity of minimal log discrepancies, extending a
previous conjecture of Shokurov. This conjecture
was proved on an ambient smooth variety in \cite{EMY}.
We apply the above theorem to 
extend this to the local complete intersection case.

\begin{theorem}\label{semicont}
If $X$ is a normal, local complete intersection variety, 
and if $Y=\sum_iq_i\cdot Y_i$ is as above, then
the function $x\longrightarrow\mld(x;X,Y)$, $x\in X$,
is lower semicontinuous.
\end{theorem}

Using Theorem~\ref{inv_adj0} and the description of minimal log
discrepancies in terms of spaces of arcs, 
we deduce also the following theorem.

\begin{theorem}\label{terminal0}
Let $X$ be a normal, local complete intersection variety.
X has log canonical (canonical, terminal) singularities if and only if
$X_m$ is equidimensional (respectively irreducible, normal) for every $m$.
\end{theorem}

The case of canonical singularities had been
conjectured by Eisenbud and Frenkel and proved
in \cite{mustata2}; the characterization
for log canonical singularities appeared there as a conjecture. 
The statement on terminal singularities answers a question
of Mirel Caib\v{a}r. All these characterizations have been proved
in \cite{EMY} in the case when $X$ is a divisor on a smooth variety.

\smallskip

The main ingredient in the proof of Theorem~\ref{inv_adj0}
is the characterization of minimal log discrepancies in terms of
jet schemes from \cite{EMY} (see Theorem~\ref{char_log_canonical} below).
In order to apply it, we need to estimate the dimension of the
subset of a jet scheme, which consists of those jets which can be lifted to
the space of arcs. This is taken care of in the next section, where
we give a small set of equations for this set. In the last section,
we use this to give the proofs of the results stated above.

\section{Equations for liftable jets}

For definitions and basic properties of jet schemes, we refer 
to \cite{mustata2}. If $X$ is a scheme of finite type over $\CC$,
then we denote its $m$th jet scheme by $X_m$, and its space of arcs
by $X_{\infty}$. When we consider $X_m$ and $X_{\infty}$, we 
restrict ourselves to the $\CC$-valued points.
If $\psi^X_m : X_{\infty}\longrightarrow X_m$
is the canonical projection, then $X_{m,\infty}$ stands
for the image of $\psi^X_m$.

We work in the ambient space $\AA^n$,
and we identify $(\AA^n)_m$ with $(\CC[t]/(t^{m+1}))^n\simeq\AA^{(m+1)n}$.
The order of $u=(u_i)_i\in(\AA^n)_m$ along an ideal 
$I\subseteq\CC[T_1,\ldots,T_n]$ is the minimum of the orders of 
$P(u)\in\CC[t]/(t^{m+1})$, over all $P\in I$. Note that either this order
is no larger than $m$, or it is equal to infinity. 

As a motivation for what follows, we mention that for
the proof of Theorem~\ref{inv_adj0} we will need to compare the dimension
of a suitable subset of $X_{m,\infty}$ with the dimension
of its intersection with $D_{m,\infty}$.  To get the right estimate, it would
be enough to cut $D_{m,\infty}$ in $X_{m,\infty}$ by a small number of equations.
It turns out to be easier
to embed (locally) $X$ in an affine space $\AA^n$, and to 
get a small number of equations for $D_{m,\infty}$ in $(\AA^n)_{m}$.
As we will see, this will be enough for our purpose.

Returning to our setting, let 
$D=V(F_1,\ldots,F_r)\subset\AA^n$, where $r\leq n$, and 
$F_i\in \CC[T_1,\ldots,T_n]$, for all $i$. We will denote
by $F$ the $r\times 1$ matrix with entries $F_i$.
Fix $\be\in\NN^r$ such that $\be_r\geq\ldots\geq\be_1$,
and fix also $p\geq\be_r$.

Up to a reordering of the variables, it is enough to consider
the following situation.
Let $(\AA^n)^{(\be)}_p$ be the set of jets in $(\AA^n)_p$ 
satisfying the following conditions.
The order along the ideal of $r$-minors of $(\partial F_k/\partial T_l)_{k\leq r,
l\leq n}$ is equal to the order along $\det(\partial F_k/\partial T_l)_{k,l\leq r}$,
and it is equal to $\be_r$; for all $i\leq r-1$, the order along the ideal
of $i$-minors of $(\partial F_k/\partial T_l)_{k\leq i,l\leq i+1}$
is equal to the order along $\det(\partial F_k/\partial T_l)_{k,l\leq i}$, and it
is equal to $\be_i$.
It is clear that $(\AA^n)_p^{(\be)}$ is a locally closed
subset of $(\AA^n)_p$.
We put $D_p^{(\be)}:=D_p\cap (\AA^n)_p^{(\be)}$, and let
$D^{(\be)}_{p,\infty}=D^{(\be)}_p\cap D_{p,\infty}$.
Our goal in this setion is to give equations for $D^{(\be)}_{p,\infty}$ 
inside $(\AA^n)_p^{(\be)}$.
Of course, we have $r(p+1)$ equations which cut
out $D_p$ in $(\AA^n)_p$,
and we will show that in addition to these, we need only $\be_r$ extra equations.

We fix first some notation. Let $u\in (\CC[t]/(t^{p+1}))^n$ be a point
in $(\AA^n)_p^{(\be)}$. We denote by $\tilde{u}\in (\CC[[t]])^n$ the unique
lifting of $u$ to polynomials of degree no more than $p$.
If we denote by $M$ the $r\times r$ matrix which is the classical adjoint of
$(\partial F_k/\partial T_l)_{k,l\leq r}$, then
we have the following

\begin{lemma}\label{equations11}{\rm (\cite{denef})}
With the above notation, $u\in (\AA^n)_p^{(\be)}$ is in
$D^{(\be)}_{p,\infty}$ if and only if
\begin{equation}\label{condition1}
{\rm ord}(M(\tilde{u})F(\tilde{u}))\geq p+\be_r+1.
\end{equation}
\end{lemma}

\begin{proof}
 It is clear that if ${\rm Jac}_F$ is the Jacobian matrix of $F$, then
$$M(\tilde{u})\cdot 
{\rm Jac}_F(\tilde{u})=(\det(\partial F_k/\partial T_l)_{k,l\leq r}
(\tilde{u})\cdot I_r, B(\tilde{u})),$$
for some $r\times (n-r)$ matrix $B$.
An easy computation based on Cramer's rule shows that
all the entries 
of $B(\tilde{u})$ have order at least $\be_r$.

We have $u\in D^{(\be)}_{p,\infty}$ if and only if
there is $v\in (t\CC[[t]])^n$, such that $F(\tilde{u}+t^pv)=0$.
Expanding using Taylor's formula gives
$$F(\tilde{u}+t^pv)=F(\tilde{u})+t^p{\rm Jac}_F(\tilde{u})v+t^{2p}(\ldots),$$
where all the entries in $(\ldots)$ have order at least $2$.
As ${\rm det}(M(\tilde{u}))\neq 0$, we have $F(\tilde{u}+t^pv)=0$
if and only if
$M(\tilde{u})\cdot F(\tilde{u}+t^pv)=0$.

Since $2p+2\geq p+\be_r+2$, and since $v\equiv 0 ({\rm mod}\,t)$, and
$M(\tilde{u}){\rm Jac}_F(\tilde{u})\equiv 0 ({\rm mod}\,t^{\be_r})$, we deduce that
if there is $v$ as above with $F(\tilde{u}+t^pv)=0$, then 
$$\ord(M(\tilde{u})\cdot F(\tilde{u}))\geq p+\be_r+1.$$

To see that this is also sufficient, note that 
$M(\tilde{u})\cdot F(\tilde{u}+t^pv)=0$
can be rewritten as $H_1(v)=\ldots=H_r(v)=0$, for
suitable polynomials $H_1,\ldots,H_r$ with coefficients in
$\CC[[t]]$.
As $H(0)\equiv 0({\rm mod}\,t)$, and
the order of an $r$-minor of the Jacobian of $H$ is zero, the existence of $v$
follows from (\ref{condition1}) via
the Implicit Function Theorem for formal power series.
\end{proof}

From now on we will assume that $u\in D_p$, 
so that $\ord(F(\tilde{u}))\geq p+1$.
Therefore the condition in (\ref{condition1}) amounts to the
vanishing of the coefficients of $t^i$ in the corresponding $r$
equations, for $p+1\leq i\leq p+1+\be_r$. Therefore we get $r\be_r$
equations. Note that while (\ref{condition1}) is written in terms of
$\tilde{u}$, since the lifting $\tilde{u}$ of $u$ was canonical, 
the resulting equations lie in the coordinate ring of $(\AA^n)_p$. 
In the remainder of this section we will show  
that, in fact, we may replace the above $r\be_r$ equations
by just $\be_r$ equations.
This is the content of the following

\begin{theorem}\label{equations12}
$D_{p,\infty}^{(\be)}$ is cut out in $(\AA^n)_p^{(\be)}\cap D_p$
by $\be_r$ equations.
\end{theorem}

\begin{proof}
We prove  the theorem by induction on $r$, the case $r=1$ being 
an immediate consequence of Lemma~\ref{equations11}.
Let $J$ denote the matrix $(\partial F_k/\partial T_l)_{k,l\leq r}$.
We denote by $\delta_{i,j}$ the $(r-1)$-minor of $J$, obtained
by deleting the $i$th row and the $j$th column. Similarly,
the $(r-2)$-minor obtained by deleting the $i$th and the $k$th rows
and the $j$th and the $l$th columns will be denoted by $\delta_{ik,jl}$.
Here we assume $i\neq k$ and $j\neq l$. We will use the following lemma,
whose proof will be given at the end of this setion.

\begin{lemma}\label{minors}
If $A$ is an arbitrary $r\times r$ matrix, and 
if $i<k$ and $j<l$, then
$$\delta_{i,j}\cdot\delta_{k,l}-\delta_{i,l}\cdot\delta_{k,j}={\rm det}(A)\cdot
\delta_{ik,jl},$$
where the $\delta$'s denote the corresponding minors of $A$, as above.
\end{lemma}

Suppose that $u\in (\AA^n)_p^{(\be)}\cap D_p$.
Recall that $M$ is the classical adjoint of $J$,
so $m_{ij}=(-1)^{i+j}\delta_{j,i}$.
The condition in (\ref{condition1}) can be written as a set of
$r$ order conditions, the $r$th one being
\begin{equation}\label{cond3}
\ord\left((-1)^r\delta_{r,r}(\tilde{u})F_r(\tilde{u})+
\sum_{j=1}^{r-1}(-1)^j\delta_{j,r}(\tilde{u})F_j(\tilde{u})\right)\geq p+\be_r+1.
\end{equation}

For every $i\leq r-1$, after multiplying the corresponding expression
from (\ref{condition1}) by
$\delta_{r,r}(\tilde{u})$ (whose order is $\be_{r-1}$), we can rewrite the 
$i$th condition in (\ref{condition1}) as
\begin{equation}\label{cond_i}
\ord\left((-1)^r\delta_{r,r}(\tilde{u})\delta_{r,i}(\tilde{u})F_r
(\tilde{u})+
\delta_{r,r}(\tilde{u})\cdot\sum_{j=1}^{r-1}(-1)^j\delta_{j,i}(\tilde{u})
F_j(\tilde{u})\right)\geq p+\be_r+\be_{r-1}+1.
\end{equation}

Granting (\ref{cond3}), which we multiply by $\delta_{r,i}(\tilde{u})$
(whose order is at least $\be_{r-1}$), we see that the condition in 
$(\ref{cond_i})$ is equivalent with
$$
\ord\left(\sum_{j=1}^{r-1}(-1)^j\delta_{j,i}(\tilde{u})\delta_{r,r}(\tilde{u})
F_j(\tilde{u})+\sum_{j=1}^{r-1}(-1)^{j+1}
\delta_{r,i}(\tilde{u})\delta_{j,r}(\tilde{u})
F_j(\tilde{u})\right)\geq p+\be_r+\be_{r-1}+1.
$$

Using Lemma~\ref{minors}, this can be further rewritten as
\begin{equation}\label{cond_iii}
\ord\left(\sum_{j=1}^{r-1}(-1)^j\delta_{jr,ir}(\tilde{u})F_j(\tilde{u})\right)
\geq p+\be_{r-1}+1.
\end{equation}
Here we used the fact that $\ord({\rm det}(J))=\be_r$.
If we denote by $J'$ the matrix formed by the first $(r-1)$
rows and columns of $J$, and by $M'$ the classical adjoint of $J'$,
and if we put $F'=(F_1,\ldots,F_{r-1})$,
then the above conditions (for $i\leq r-1$) can be rewritten as
$$\ord(M'(\tilde{u})\cdot F'(\tilde{u}))\geq p+\be_{r-1}+1.$$
By induction on $r$, this is equivalent with a set of $\be_{r-1}$ equations
$G_1({u})=\ldots=G_{\be_{r-1}}({u})=0$. 

Given this, we go back to the condition in (\ref{cond3}).
Since $u\in D_p$, we have
${\rm ord}\,F(\tilde{u})\geq p+1$, so that 
$\ord(\delta_{r,r}(\tilde{u})F_r(\tilde{u}))\geq p+\be_{r-1}+1$.
Moreover, by expanding each $(r-1)$-minor
along the last row, we write
$$\sum_{j=1}^{r-1}(-1)^j\delta_{j,r}(\tilde{u})F_j(\tilde{u})
=\sum_{i=1}^{r-1}(-1)^{i+r+1}\frac{\partial F_r}{\partial T_i}(\tilde{u})
\cdot\sum_{j=1}^{r-1}
(-1)^j\delta_{jr,ir}(\tilde{u})F_j(\tilde{u}),$$
which by (\ref{cond_iii}) has order at least $p+\be_{r-1}+1$. Therefore, in order
to ensure $(\ref{cond3})$ we need only $(\be_r-\be_{r-1})$ 
more equations $G_{\be_{r-1}+1}(u)=\ldots=G_{\be_r}(u)=0$. 
This concludes the proof
of the theorem.
\end{proof}

We give now the proof of the lemma we have used.

\begin{proof}[Proof of Lemma~\ref{minors}]
All the determinants are considered as
polynomial functions on the affine space of all matrices.
Let $P$ denote the difference of products of minors in the left hand side
of the formula in the lemma. To see that $P(A)$ is divisible
by $\det(A)$ for every matrix $A$, it is enough to show that if
$A$ is such that one column is a linear 
combination of the other columns, then $P(A)=0$. Let that column be indexed 
by $m$, so we can write $C_m=\sum_{\alpha\neq m}\lambda_{\alpha}C_{\alpha}$
($C_{\alpha}$ denotes the $\alpha$th column). 
If $m\neq j$, $l$, then $\delta_{i,j}(A)=(-1)^{m-j}\lambda_j\delta_{i,m}(A)$,
$\delta_{k,l}(A)=(-1)^{m-l}\lambda_l\delta_{k,m}(A)$,
$\delta_{i,l}(A)=(-1)^{m-l}\lambda_l\delta_{i,m}(A)$,
and $\delta_{k,j}(A)=(-1)^{m-j}\lambda_j\delta_{k,m}(A)$. This 
clearly gives $P(A)=0$.

Suppose now, for example, that $m=j$. We have $\delta_{k,l}(A)=
(-1)^{m-l}\lambda_l\delta_{k,j}(A)$ and $\delta_{i,l}(A)
=(-1)^{m-l}\lambda_l\delta_{i,j}(A)$, and this gives $P(A)=0$.
The case $m=l$ is similar.

It follows from the above discussion that
we can write $P(A)={\rm det}(A)\cdot Q(A)$. Since the degree
of $P$ with respect to any of $a_{i\beta}$, $a_{\beta j}$,
$a_{k\beta}$ and $a_{\beta l}$ is at most one, for all $1\leq\beta\leq r$,
we see that $Q$ depends only on the remaining variables.

Therefore, if $B$ is the matrix obtained from a matrix $A$ by deleting
the $i$th and the $k$th rows and the $j$th and the $l$th columns, then 
$Q$ depends only on $B$. In order to compute $Q(B)$, we may assume that
$a_{ij}=a_{kl}=1$, and that all the other entries on those two rows 
and columns are zero. In this case, it is clear that
$\delta_{i,j}(A)=(-1)^{k+l}\delta_{ik,jl}(A)$,
$\delta_{k,l}(A)=(-1)^{i+j}\delta_{ik,jl}(A)$,
while $\delta_{i,l}(A)=0=\delta_{k,j}(A)$. Since
${\rm det}(A)=(-1)^{k+l+i+j}\delta_{ik,jl}(A)$, we see that
$Q(A)=\delta_{ik,jl}(A)$, as required.
\end{proof}

\begin{remark}
Note that the above description of the equations
does not assume that $D$ is locally 
complete intersection, but only that the number
of equations is not greater than the dimension of the affine space.
\end{remark}

\section{Log discrepancies and Inversion of Adjunction}

We start by recalling the definition of the version of minimal log
discrepancies from \cite{EMY}.
Consider a pair $(X,Y)$, where $X$ is a $\QQ$-Gorenstein normal variety, and
where $Y$ stands for a formal combination
$\sum_{i=1}^kq_i\cdot Y_i$, where $q_i\in\RR_+$
and $Y_i\subset X$ are proper closed subschemes. 
A divisor $E$ over $X$ is a prime divisor on some smooth variety $X'$,
such that we have a proper, birational morphism $f : X'\longrightarrow X$.
In fact, we identify two such divisors if they give the same valuation
of the function field of $X$. In particular, we may assume that the inverse image
of each $Y_i$ on $X'$ is an effective Cartier divisor.
The center of $E$ on $X$ is $c_X(E):=f(E)$.

The log discrepancy $a(E;X,Y)$ is defined such that $a(E;X,Y)-1$
is the coefficient of $E$ in $K_{X'/X}-f^{-1}(Y)$. Here 
$f^{-1}(Y)=\sum_iq_if^{-1}(Y_i)$, where $f^{-1}(Y_i)$ is the scheme
theoretic inverse image of $Y_i$. Recall the convention that $K_{X'/X}$
is a divisor supported on the exceptional locus of $f$.

If $W\subseteq X$ is a closed subset, then the minimal log discrepancy
of $(X,Y)$ on $W$ is defined by
$$\mld(W;X,Y):=\min_{c_X(E)\subseteq W}a(E;X,Y).$$
For more on these invariants, we refer to
\cite{EMY}, Section 1, or to \cite{ambro}.

We mention that the pair $(X,Y)$ is called log canonical if
$\mld(X;X,Y)\geq 0$. This is a requirement for having such interesting
invariants: if $\dim\,X\geq 2$, and if $(X,Y)$ is not log canonical,
then $\mld(X;X,Y)=-\infty$. Note also that if $\mld(W;X,Y)\geq 0$,
then the pair $(X,Y)$ is log canonical in some neighbourhood of $W$.

We need some preparations in order to state the characterization of
minimal log discrepancies from \cite{EMY}.
We assume from now on that $X$ is a normal, locally complete
intersection variety, as this is the setting we have in Theorem~\ref{inv_adj0}.
Let $d=\dim\,X$, and
let $Z\subset X$ be the Jacobian subscheme of $X$, i.e., the subscheme
defined by the $d$th Fitting ideal of $\Omega_X^1$. Note that the
support of $Z$ is the singular locus $X_{\rm sing}$ of $X$.

If $Y\subset X$ is a closed subscheme, then we have a function
$F_Y : X_{\infty}\longrightarrow\NN\cup\{\infty\}$, such that
$F_Y(\gamma)$ is 
the order of vanishing of $\gamma$ along the ideal of $Y$.
Recall that we have canonical projections $\psi^X_m : X_{\infty}\longrightarrow 
X_m$, so that $F_Y^{-1}(\geq m)=(\psi_{m-1}^X)^{-1}(Y_{m-1})$, for all $m\geq 1$.
 
It is a theorem of Denef and Loeser from \cite{denef} that if $m\geq e$,
then the induced projections
$$\psi_{m+1}^X(F_Z^{-1}(e))\longrightarrow\psi_m^X(F_Z^{-1}(e))$$
are locally trivial, with fiber $\AA^d$. Suppose now that
we have proper closed subschemes $Y_1,\ldots,Y_k\subset X$,
and suppose that $A=\bigcap_iF_{Y_i}^{-1}(\geq m_i)$, for some $m_i\in\NN$.
The codimension of $A$ in $X_{\infty}$ is defined by
\begin{equation}\label{def_codim}
\codim(A, X_{\infty}):=\min_{e,m}\{(m+1)d-\dim\,\psi_m^X
(A\cap F_Z^{-1}(e))\},
\end{equation}
where the minimum is over all $e\in\NN$ and $m\geq \max\{e,m_1,\ldots,m_k\}$ 
(the above result of
Denef and Loeser implies that 
the  expression in (\ref{def_codim}) is independent of which
$m$ 
we choose, as above). We use the convention $\dim(\emptyset)=-\infty$.

\begin{theorem}\label{char_log_canonical}{\rm (\cite{EMY})}
Let $(X,Y)$ be a pair, with $X$  a normal, locally complete intersection
variety, and $Y=\sum_{i=1}^kq_i\cdot Y_i$, where $q_i\in\RR_+$, and where
$Y_i\subset X$ are proper closed subschemes. 
Fix a proper closed subscheme $W\subset X$ and let $\tau\in\RR_+$.
If $Z\subset X$ is the Jacobian
subscheme, then
$\mld(W;X,Y)\geq\tau$ if and only if
$$\codim\left(F_W^{-1}(\geq 1)
\cap F_Z^{-1}(\geq e)\cap\bigcap_{i=1}^k F_{Y_i}^{-1}(\geq m_i), X_{\infty}\right)
\geq e+\sum_{i=1}^kq_im_i+\tau,$$
for every $e\in\NN$, $(m_i)_i\in\NN^k$.
\end{theorem}

\bigskip

We combine now the above theorem with our
result from the previous section to give
a proof of the case of the  Inversion of Adjunction Conjecture
which was stated in Section 1.

\begin{proof}[Proof of Theorem~\ref{inv_adj0}]
The inequality 
\begin{equation}\label{adjunction}
\mld(W;X,D+Y)\leq\mld(W;D,Y\vert_D)
\end{equation}
is well-known in general and follows by adjunction (see, for example,
the proof of Proposition~7.3.2 in \cite{kollar}, or the proof
of Theorem 1.6 in \cite{EMY}). 
We  give the proof of the reverse inequality.

Working locally, we may assume that $X$ is a subvariety of codimension $r-1$ of
an open subset $U$ of an affine space $\AA^n$. Moreover, we may assume that
$X=H_1\cap\ldots H_{r-1}$, and $D=X\cap H_r$, where $H_i$ is defined in $U$
by $F_i\in\CC[T_1,\ldots,T_n]$. An inductive application of
(\ref{adjunction}) gives
$$\mld(W;U,Y+\sum_{i=1}^r H_i)\leq\mld(W;X,D+Y).$$ Therefore it is enough
to prove $\mld(W;D,Y\vert_D)\leq\mld(W;U,Y+\sum_{i=1}^r H_i)$. 

Consider $\tau\in\RR_+$, and suppose that
$\mld(W; D, Y\vert_D)\geq\tau$. We assume that $\mld(W;U,Y+\sum_{i=1}^r H_i)
<\tau$, and we will derive a contradiction. We first apply 
Theorem~\ref{char_log_canonical} to the smooth variety $U$, to deduce that
we can find $m_1,\ldots,m_k\in\NN$, and $m'_1,\ldots,m'_r\in\NN$
such that 
$$\codim\left(F_W^{-1}(\geq 1)\cap
\bigcap_{i=1}^kF_{Y_i}^{-1}(\geq m_i)\cap\bigcap_{j=1}^r
F_{H_j}^{-1}(\geq m'_j),U_{\infty}\right)
<\sum_{i=1}^km_iq_i+\sum_{j=1}^rm'_j+\tau.$$
Let us choose $m\geq m_i$, $m'_j$, for all $i$ and $j$,
 and consider the projection to $U_m$
of the above subset of $U_{\infty}$. We may choose an irreducible component
$V$ of this projection, such that
$\dim\,V>(m+1)n-\sum_iq_im_i-\sum_jm'_j-\tau$.

As in Lemma~3.2 from \cite{EMY}, one can show that 
$V\cap\psi_m^D(D_{\infty})\neq\emptyset$. Moreover, given
$v\in V\cap(\psi_m^D)(D_{\infty})$, there is a 
$\tilde{v}\in(\psi_m^D)^{-1}(v)\setminus Z_{\infty}$.
Here $Z\subset D$ is the Jacobian subscheme of $D$,
which is cut out on $D$ by the $r$-minors of the matrix
$(\partial F_i/\partial T_j)_{i\leq r, j\leq n}$.
If $f=F_Z(\tilde{v})$, then
we may assume that $f\leq m$. Indeed, otherwise, we may
just replace $m$ by some $m'\geq f$, $V$ by its inverse image in
$X_{m'}$, and $v$ by $\psi^D_{m'}(\tilde{v})$.  

Choose now $\be_r$ to be the smallest order of vanishing along $Z$
of an element in $V\cap\Image(\psi_m^D)$. We see that
$\be_r\leq f\leq m$. After renumbering the variables, we may
assume that some element in $(\psi_m^D)^{-1}(V)$
vanishes along  
$\det(\partial F_i/\partial X_j)_{i,j\leq r}$
with order $\be_r$. For these elements, we choose $\be_{r-1}\leq\be_r$
to be the smallest order of vanishing along the ideal generated by the
$(r-1)$-minors of $(\partial F_i/\partial X_j)_{i\leq r-1,j\leq r}$.
Moreover, after renumbering the variables we may assume that
this minimum is achieved by an element as above, along
$\det(\partial F_i/\partial X_j)_{i,j\leq r-1}$. 
We proceed in this way to obtain a sequence 
$\be_1\leq\be_2\leq\ldots\leq\be_r$.
We fix also $\tilde{v}'\in(\psi_m^D)^{-1}(V\cap(\AA^n)_m^{(\be)})$,
where we use the notation introduced in the previous section.

By Greenberg's theorem (see \cite{Gr}), there is a
$p\geq m$ such that every element in $D_m$ which can be lifted
to $D_p$ can be lifted to $D_{\infty}$. 
If $\phi : D_p\longrightarrow D_m$ is the canonical projection,
then let $V'$ be an irreducible component
of $\phi^{-1}(V\cap D_m)$ which
contains $\psi_p^D(\tilde{v}')$. As
 $\phi^{-1}(V\cap D_m)$
is cut out in $V\times\AA^{(p-m)n}$ by $\sum_{j=1}^r(p-m'_j+1)$ equations, we
see that
$$\dim\,V'> (p+1)(n-r)-\sum_{i=1}^kq_im_i-\tau.$$

It follows from our choice of $\be$ and $p$ that $V'_0:=V'\cap
(\AA^n)_p^{(\be)}$ is open in $V'$. 
Moreover, Theorem~\ref{equations12}
shows that $V''_0:=V'_0\cap\Image(\psi_p^D)$ is cut out in
$V'_0$ by $\be_r$ equations. Since $V''_0$ is nonempty, by construction,
we see that 
$$\codim((\psi_p^D)^{-1}(V''_0),D_{\infty})<\be_r+\sum_{i=1}^kq_im_i+\tau,$$
while $(\psi_p^D)^{-1}(V''_0)\subseteq F_W^{-1}(\geq 1)\cap F_Z^{-1}(\geq \be_r)
\cap\bigcap_iF_{Y_i}^{-1}(\geq
m_i)$. This contradicts $\mld(W;D, Y\vert_D)\geq\tau$, by 
Theorem~\ref{char_log_canonical}, and completes the proof.

\end{proof}

\begin{corollary}\label{invadj_lci}
Let $V$ be a smooth variety, $X\subset V$ a normal, local complete intersection
subvariety of codimension $r$, and $Y=\sum_iq_i\cdot Y_i$, where
$q_i\in\RR_+$ and where $Y_i\subset X$ are proper closed subschemes.
If $W\subset X$ is a proper closed subset, then
$\mld(W;X,Y)=\mld(W;V,Y+r\cdot X)$.
\end{corollary}

\begin{proof}
Working locally, we may assume that $X=H_1\cap\ldots\cap H_r$, where
the $H_i$'s are effective divisors on $V$. Note that $H_1\cap\ldots\cap H_{i}$
is normal around $X$, for every $i$.
It is enough
to show that we have the inequalities
$$\mld(W;V,Y+\sum_{i=1}^rH_i)\leq\mld(W;V,Y+r\cdot X)\leq\mld(W;X,Y),$$
as the extremal terms are equal by an inductive application of
Theorem~\ref{inv_adj0}. The first inequality follows immediately
from the definition of minimal log discrepancies.

For the second one, we proceed as follows. Let $\pi : \widetilde{V}
\longrightarrow V$ be the blowing-up of $V$ along $X$, with
exceptional divisor $E$. If $\pi' : E\longrightarrow X$ is the restriction 
of $\pi$, then $\pi'$ is a projective bundle with fiber ${\mathbb P}^{r-1}$,
so that $\mld(W;X,Y)=\mld(\pi^{-1}(W);E,\pi^{-1}(Y))$.
Moreover, $E$ is reduced and irreducible. It is also locally a complete intersetion,
hence so is $\widetilde{V}$. By computing $K_{\widetilde{V}/V}$
at the general point of $E$, we see that $K_{\widetilde{V}/V}=(r-1)E$.

From the formula describing the behaviour of minimal log
discrepancies under a proper, birational morphism
(see \cite{EMY}, Proposition~1.3 (iv)),
we get
$$\mld(W;V,Y+r\cdot X)=\mld(\pi^{-1}(W);\widetilde{V},
\pi^{-1}(Y)+rE-K_{\widetilde{V}/V})$$
$$=\mld(\pi^{-1}(W);\widetilde{V},
\pi^{-1}(Y)+E)\leq
\mld(\pi^{-1}(W);E,\pi^{-1}(Y)),$$
where the inequality follows from (\ref{adjunction}).
Therefore $\mld(W;V,Y+r\cdot X)\leq\mld(W;X,Y)$, so we are done.
\end{proof}

Together with the corresponding result from \cite{EMY} for the smooth case,
this immediately gives the semicontinuity of minimal log discrepancies
for local complete intersetion varieties.

\begin{proof}[Proof of Theorem~\ref{semicont}]
Recall that $X$ is a normal, local complete intersection variety.
Working locally, we may assume that $X$ is a codimension $r$
subvariety of a smooth variety $V$. It follows from the above corollary
that for every $x\in X$, we have
$$\mld(x;X,Y)=\mld(x;V,Y+r\cdot X),$$
so it is enough to use the corresponding 
semicontinuity assertion on smooth varieties.
This is Theorem~0.3 in \cite{EMY}.
\end{proof}

We give now the proof of the characterization 
of certain classes of singularities for local complete intersection
varieties.

\begin{proof}[Proof of Theorem~\ref{terminal0}]

Recall the definition of log canonical, canonical, and
terminal singularities in terms of minimal log discrepancies.
A normal,
$\QQ$-Gorenstein 
variety $X$ has log canonical (canonical, terminal) singularities
if and only if we have $\mld(X_{\rm sing}; X, \emptyset)\geq 0$
(respectively, $\geq 1$, $>1$). Note also that if $X$ is
Gorenstein, as in our case, then the above minimal log discrepancy
is an integer.

We consider first the characterization of log canonical singularities.
Let $\pi_m : X_m\longrightarrow X$ be the canonical projection.
$X_m$ is equidimensional if and only if $\dim\pi_m^{-1}(X_{\rm sing})
\leq (m+1)\dim\,X$. Moreover, if this is the case, then $X_m$
is locally a complete intersection and $\dim\,X_m=(m+1)\dim\,X$
(see \cite{mustata2}). 
Working locally, we may assume that $X$ is a codimension $r$
subvariety of a smooth variety $V$. It follows from 
Theorem~\ref{char_log_canonical} that $X_m$ is equidimensional 
for every $m$ if and only if
$\mld(X_{\rm sing};V, r\cdot X)\geq 0$. By Corollary~\ref{invadj_lci},
this is equivalent with $\mld(X_{\rm sing};X,\emptyset)\geq 0$, i.e.,
with $X$ having log canonical singularities. This completes this case.

Similarly, $X_m$ is irreducible if and only if $\dim\pi_m^{-1}(X_{\rm sing})
\leq (m+1)\dim\,X-1$ (see \cite{mustata2}). 
The characterization of canonical singularities follows as above.
From now on, we may assume that $X$ has canonical singularities.
In particular, $X_m$ is locally a complete intersection for every $m$.
Serre's criterion implies that $X_m$ is normal if and only if
$\dim (X_m)_{\rm sing}\leq (m+1)\dim\,X-2$. One can easily show
that because $X_m$ is equidimensional, we have $(X_m)_{\rm sing}=
\pi_m^{-1}(X_{\rm sing})$ (see the proof of Theorem~3.3
in \cite{EMY}). Therefore $X_m$ is normal for every $m$
if and only if $\mld(X_{\rm sing}; V, r\cdot X)\geq 2$, and we conclude as above,
applying Corollary~\ref{invadj_lci}.
\end{proof}

\begin{remark}
It is clear that if $X$ is a variety such that $X_1$ is irreducible,
then $\dim\,X_{\rm sing}\leq\dim\,X-2$. If $X$ is locally a complete
intersection,
then $X$ has to be normal. It follows that in the above characterization
for canonical or terminal singularities, we do not need to assume that
$X$ is normal.
\end{remark}

\bigskip

\subsection*{Acknowledgements}
We are grateful to Florin Ambro, Vyacheslav V.~Shokurov, and
Takehiko Yasuda for useful discussions. 
The first author was partially
supported by NSF grant DMS 0200278. The second author
served as a Clay Mathematics Institute Long-Term Prize Fellow
while this research has been done.

\providecommand{\bysame}{\leavevmode \hbox \o3em
{\hrulefill}\thinspace}


\begin{thebibliography}{Gr}


\bibitem[Am]{ambro}
F.~Ambro, On minimal log discrepancies, Math. Res. Letters \textbf{6}
(1999), 573--580.

\bibitem[DL]{denef}
J.~Denef and F.~Loeser, Germs of arcs on singular varieties
and motivic integration, Invent. Math. \textbf{135} (1999), 201--232.


\bibitem[EMY]{EMY}
L.~Ein, M.~Musta\c{t}\v{a} and T.~Yasuda,
Log discrepancies, jet schemes and Inversion of Adjunction,
Invent. Math., to appear.

\bibitem[Gr]{Gr}
M.~Greenberg, Rational points in henselian discrete valuation rings,
Publ. Math. I.H.E.S. \textbf{31} (1966), 59--64. 


\bibitem[K+]{kollar}
J.~Koll\'{a}r (with 14 coauthors), Flips and abundance for algebraic threefolds,
Ast\'{e}risque \textbf{211}, 1992.


\bibitem[Mu2]{mustata2}
M.~Musta\c{t}\v{a}, Jet schemes of locally complete intersection
canonical singularities, with an appendix by D.~Eisenbud and E.~Frenkel,
Invent. Math. \textbf{145} (2001), 397--424.


\end{thebibliography}
\end{document}